\documentclass[12pt]{amsart}

\usepackage{amsmath}
\usepackage{amssymb}
\usepackage{amscd}
\usepackage{amsthm}
\usepackage{epsfig}
\usepackage{color}
\theoremstyle{plain}

\newtheorem{thm}[equation]{Theorem}

\newtheorem{prop}[equation]{Proposition}

\newtheorem{lem}[equation]{Lemma}

\newtheorem{cor}[equation]{Corollary}

\theoremstyle{remark}

\newtheorem{rem}[equation]{Remark}
\theoremstyle{definition}

\newtheorem{convention}[equation]{Convention}

\newtheorem{defn}[equation]{Definition}

\newtheorem{question}[equation]{Question}

\newcommand{\N}{\mathbb N}

\newcommand{\R}{\mathbb R}

\newcommand{\Z}{\mathbb Z}

\newcommand{\length}{\operatorname{length}}

\newcommand{\al}{\alpha}

\def\eps{\epsilon}
\def\ga{\gamma}
\def\Ga{\Gamma}

\def\la{\lambda}
\def\La{\Lambda}

\def\<{\langle}
\def\>{\rangle}

\def\om{\omega}

\def\si{\sigma}
\def\Si{\Sigma}

\def\be{\beta}
\def\del{\delta}

\def\ulim{\mathop{\hbox{$\om$-lim}}}

\newcommand{\ol}{\overline}

\def\om{\omega}

\begin{document}

\title{Triangle inequalities in path metric spaces}

\author{Michael Kapovich}
\date{November 5, 2006}
\address{Department of Mathematics, University of California,
Davis, CA 95616} \email{kapovich@math.ucdavis.edu}

\begin{abstract}
We study side-lengths of triangles in path metric spa\-ces. We prove that unless such a space $X$ is bounded, or 
quasi-isometric to $\R_+$ or to $\R$, every triple of real numbers satisfying the strict triangle inequalities, 
is realized by the side-lengths of a triangle in $X$. We construct an example of a complete path metric space quasi-isometric 
to $\R^2$ for which every degenerate triangle has one side which is shorter than a certain uniform constant. 
\end{abstract}

\maketitle

Given a metric space $X$ define
$$
K_3(X):=\{ (a, b, c)\in \R_+^3: \exists \hbox{~~points~~} x, y, z
\hbox{~~so that~~}\\$$
$$
 d(x,y)=a, d(y,z)=b, d(z, x)=c\}.
$$
Note that $K_3(\R^2)$ is the closed convex cone $K$ in $\R_+^3$
given by the usual triangle inequalities. On the other hand,  if $X=\R$ then $K_3(X)$
is the boundary of $K$ since all triangles in $X$ are degenerate. 
If $X$ has finite diameter, $K_3(X)$ is a bounded set.

In \cite[Page 18]{Gromov00} Gromov raised the following

\begin{question}
Find {\em reasonable} conditions on path metric spaces $X$, under
which $K_3(X)=K$.
\end{question}

It is not so difficult to see that for a path metric space $X$ quasi-isometric to
$\R_+$ or $\R$, the set $K_3(X)$ does not contain the interior of $K$, 
see Section \ref{exceptional}. Moreover, every triangle in such $X$ is $D$-degenerate for 
some $D<\infty$ and therefore $K_3(X)$ is contained in the $D$-neighborhood of $\partial K$.

Our main result  is essentially the converse to the above observation:  

\begin{thm}\label{main}
Suppose that $X$ is an unbounded path metric space 
not quasi-isometric to $\R_+$ or $\R$. Then:

1. $K_3(X)$ contains the interior of the cone $K$.

2. If, in addition, $X$ contains arbitrary long geodesic segments, 
then $K_3(X)= K$.
\end{thm}

In particular, we obtain a complete answer to Gromov's question for
geodesic metric spaces, since an unbounded geodesic metric space
clearly contains arbitrarily long geodesic segments. In Section
\ref{examples}, we give an example of a (complete) path metric
space $X$ quasi-isometric to $\R^2$, for which
$$
K_3(X)\ne K. 
$$
Therefore, Theorem \ref{main} is the optimal result.

The proof of Theorem \ref{main} is easier under the assumption that
$X$ is a proper metric space: In this case $X$ is necessarily
complete, geodesic metric space. Moreover, every unbounded sequence
of geodesic segments $\ol{o x_i}$ in $X$ yields a geodesic ray. The
reader who does not care about the general path metric spaces can
therefore assume that $X$ is proper. The arguments using the
ultralimits are then replaced by the Arcela--Ascoli theorem.

\medskip
Below is a sketch of the proof of Theorem \ref{main} under the extra assumption that $X$ is proper. 
Since the second assertion of Theorem \ref{main} is clear, we have to  prove only the first 
statement. We define {\em $R$-tripods} $T\subset X$, as unions $\ga\cup \mu$ of two geodesic segments 
$\ga, \mu\subset X$, having the lengths $\ge R$ and $\ge 2R$ respectively, 
so that:

1. $\ga\cap \mu=o$ is the end-point of $\ga$.

2. $o$ is  distance $\ge R$ from the ends of $\mu$.

3. $o$ is a nearest-point projection of $\ga$ to $\mu$. 

The space $X$ is called $R$-{\em compressed} if it contains no $R$-tripods. 
The space $X$ is called {\em uncompressed} if it is not $R$-compressed for any $R<\infty$.

We break the proof of Theorem \ref{main} in two parts:

{\bf Part 1.} 

\begin{thm}\label{A}
If $X$ is uncompressed then $K_3(X)$ contains the interior of
$K_3(\R^2)$.
\end{thm}

The proof of this theorem is mostly topological. The side-lengths of triangles in $X$ 
determine a continuous  map
$$
L: X^3\to K
$$
Then $K_3(X)=L(X^3)$. Given a point $\kappa$ in the interior of $K$, we consider an $R$-tripod 
$T\subset X$ for sufficiently large $R$. We then restrict to triangles in $X$ with vertices in $T$. 
We construct a 2-cycle $\Si\in Z_2(T, \Z_2)$ whose image under $L_*$ 
determines a nontrivial element of $H_2(K\setminus \kappa, \Z_2)$. Since $T^3$  is contractible, 
there exists a 3-chain $\Ga\in C_3(T^3, \Z_2)$ with the boundary $\Si$. Therefore 
the support of $L_*(\Ga)$ contains the point $\kappa$, which implies that $\kappa$ belongs to 
the image of $L$. 

\begin{rem}
Gromov observed in \cite{Gromov00} that {\em uniformly contractible} metric spaces $X$ have {\em large} $K_3(X)$. 
Although uniform contractibility is not relevant to our proof, the key argument here indeed has the 
coarse topology flavor. 
\end{rem}

{\bf Part 2.} 

\begin{thm}\label{B}
If $X$ is a compressed unbounded path metric space, then $X$ is quasi-isometric to $\R$ or $\R_+$.  
\end{thm}

Assuming that $X$ is compressed, unbounded and is not quasi-iso\-met\-ric to $\R$ and to $\R_+$, 
we construct three diverging geodesic rays $\rho_i$ in $X$, $i=1, 2, 3$. Define 
$\mu_i\subset X$ to be the geodesic segment connecting $\rho_1(i)$ and $\rho_2(i)$.  
Take $\ga_i$ to be the shortest segment connecting $\rho_3(i)$ to $\mu_i$. Then $\ga_i\cup \mu_i$ 
is an $R_i$-tripod with $\lim_i R_i=\infty$, which contradicts the assumption that $X$ is compressed.

\medskip
{\bf Acknowledgements.} During this work the author was
partially supported by the NSF grants DMS-04-05180 and DMS-05-54349. 
Most of this work was done when the author 
was visiting the Max Plank Institute for Mathematics in Bonn.

\section{Preliminaries}
\label{basics}

\begin{convention}\label{conv}
 All homology will be taken with the $\Z_2$-coefficients.
\end{convention}

In the paper we will talk about {\em ends of a metric space} $X$. 
Instead of looking at the noncompact complementary components of {\em relatively compact open subsets} of $X$ 
as it is usually done for topological spaces, we will define ends of $X$ by considering 
unbounded complementary components of bounded subsets of $X$. With this modification, 
the usual definition goes through. 
 
If $x, y$ are points in a topological space $X$, we let $P(x,y)$ denote the set of 
continuous paths in $X$ connecting $x$ to $y$. For $\al\in P(x,y), \be\in P(y, z)$ we 
let $\al*\be\in P(x, z)$ denote the concatenation of $\al$ and $\be$. Given a path 
$\al: [0, a]\to X$ we let $\bar\al$ denote the reverse path
$$
\bar\al(t)= \al(a-t). 
$$

\subsection{Triangles and their side-lengths}

We set $K:=K_3(\R^2)$; it is the cone in $\R^3$ given by
$$
\{ (a, b, c): a\le b+c, b\le a+c, c\le a+b\}.
$$
We metrize $K$ by using the maximum-norm on $\R^3$. 

By a {\em triangle} in a metric space $X$ we will mean an ordered triple $\Delta=(x, y, z)\in X^3$.
We will refer to the numbers $d(x, y), d(y, z), d(z, x)$ as the {\em side-lengths} of $\Delta$, 
even though these points are not necessarily connected by geodesic segments. 
The sum of the side-lengths of $\Delta$ will be called the {\em perimeter} of $\Delta$.

We have the continuous map
$$
L: X^3\to K
$$
which sends the triple $(x, y, z)$ of points in $X$ to the triple of side-lengths
$$
(d(x, y), d(y, z), d(z, x)).
$$
Then $K_3(X)$ is the image of $L$.

Let $\eps\ge 0$. We say that a triple $(a, b, c)\in K$ is
$\eps$-{\em degenerate} if, after reordering if necessary the
coordinates $a, b, c$, we obtain
$$
a+\eps \ge b+c.
$$
Therefore every $\eps$-degenerate triple is within distance $\le
\eps$ from the boundary of $K$. A triple which is not
$\eps$-degenerate is called $\eps$-nondege\-ne\-rate.
A triangle in a metric space $X$ whose side-lengths form an $\eps$-{\em degenerate} triple,
is called $\eps$-{\em degenerate}. A $0$-degenerate triangle is called {\em degenerate}. 

\subsection{Basic notions of metric geometry}

For a subset $E$ in a metric space $X$ and $R<\infty$ we let $N_R(E)$ denote the metric $R$-neighborhood of $E$ in $X$:
$$
N_R(E)=\{ x\in X: d(x, E)\le R\}.
$$

\begin{defn}
Given a subset $E$ in a metric space $X$ and $\eps>0$, we define the
$\eps$-{\em nearest--point projection} $p=p_{E,\eps}$ as the map which
sends $X$ to the set $2^E$ of subsets in $E$:
$$
y\in p(x) \iff d(x,y)\le d(x,z)+\eps, \quad \forall z\in E.
$$
If $\eps=0$, we will abbreviate $p_{E,0}$ to $p_E$.  
\end{defn}

{\bf Quasi-isometries.} Let $X, Y$ be metric spaces. A map
$f: X\to Y$ is called an $(L,A)$-quasi-isometric embedding (for $L\ge 1$ and $A \in \R$) 
if for every pair of points $x_1, x_2\in X$ we have
$$
L^{-1} d(x_1, x_2)- A\le d(f(x_1), f(x_2))\le L d(x_1, x_2)+A.
$$

A map $f$ is called an $(L,A)$-{\em quasi-isometry}  if it is an $(L,A)$-quasi-isometric embedding so that
$N_A(f(X))=Y$. Given an $(L,A)$-quasi-isometry, we have the {\em quasi-inverse} map 
$$
\bar{f}: Y\to X
$$
which is defined by choosing for each $y\in Y$ a point $x\in X$ so that $d(f(x),y)\le A$.
The quasi-inverse map $\bar{f}$ is an $(L, 3A)$-quasi-isometry. An $(L,A)$-{\em quasi-isometric embedding} $f$ 
of an interval $I\subset \R$ into a metric space $X$ is called an $(L,A)$-{\em quasi-geodesic} in $X$. 
If $I=\R$, then $f$ is called a {\em complete} quasi-geodesic.

A map $f: X\to Y$ is called a {\em quasi-isometric embedding} (resp. a 
{\em quasi-isometry}) if it is an $(L,A)$-quasi-isometric embedding (resp.
$(L,A)$-quasi-isometry) for some $L\ge 1, A \in \R$.

Every quasi-isometric embedding $\R^n\to \R^n$ is a quasi-isometry, see for 
instance \cite{Kapovich-Leeb(1995)}.

\bigskip 
{\bf Geodesics and pathe metric spaces.} 

A {\em geodesic} in a metric space is an isometric embedding of an
interval into $X$. By abusing the notation, we will identify
geodesics and their images. A metric space is called {\em geodesic}
if any two points in $X$ can be connected by a geodesic. By abusing the notation we let $\ol{xy}$ denote 
a geodesic connecting $x$ to $y$, even though this geodesic is not necessarily unique. 

The length of a continuous curve $\ga: [a,b]\to X$ in a metric space, is defined  as
$$
\length(\ga)= \sup \Big{ \{ }  \sum_{i=1}^n d(\ga(t_{i-1}), \ga(t_{i})): a=t_0< t_1<... < t_n=b \Big{ \} }.
$$
A path $\ga$ is called {\em rectifiable} if $\length(\ga)<\infty$. 

A metric space $X$ is called a {\em path metric space} if for every
pair of points $x, y\in X$ we have
$$
d(x, y)=\inf \{ \length(\ga): \ga\in P(x, y)\}. 
$$
We say that a curve $\ga: [a, b]\to X$ is $\eps$-geodesic if
$$
\length(\ga)\le d(\ga(a), \ga(b))+\eps.
$$
It follows that every $\eps$-geodesic is $(1,\eps)$--quasi-geodesic.
We refer the reader to \cite{Gromov00, Burago-Ivanov} for the further details on path metric spaces.

\subsection{Ultralimits}

Our discussion of ultralimits of sequences of metric space will be somewhat brief, we refer the reader to 
\cite{Kapovich-Leeb(1995), Gromov00, Kapovich2000, Burago-Ivanov, Roe(2003)} for the detailed definitions and discussion. 

Choose  a nonprincipal ultrafilter $\om$ on $\N$. Suppose that we are given a sequence of pointed 
metric spaces $(X_i, o_i)$, where $o_i\in X_i$. The {\em ultralimit}
$$
(X_\om, o_\om)=\ulim (X_i, o_i)
$$
is a pointed metric space whose elements are equivalence classes of sequences $x_i\in X_i$. 
The distance in $X_\om$ is the $\om$-limit: 
$$
\ulim d(x_i, y_i). 
$$
One of the key properties of ultralimits which we will use repeatedly is the following. Suppose that $(Y_i,p_i)$ is a 
sequence of pointed metric spaces. Assume that we are  
given a sequence of $(L_i, A_i)$-quasi-isometric embeddings
$$
f_i: X_i \to Y_i
$$
so that $\ulim d(f(o)i), p_i))<\infty$ and
$$
\ulim L_i=L<\infty, \quad \ulim A_i=0. 
$$
Then there exists the ultralimit $f_\om$ of the maps $f_i$, which is an $(L,0)$-quasi-isometric embedding
$$
f_\om: X_\om\to Y_\om. 
$$
In particular, if $L=1$, then $f_\om$ is an isometric embedding.

\medskip
{\bf Ultralimits of constant sequences of metric spaces.} Suppose that $X$ is a path metric space. 
Consider the constant sequence $X_i=X$ for all $i$. 
If $X$ is a proper metric space and $o_i$ is a bounded sequence, the
ultralimit $X_\om$ is nothing but $X$ itself. In general, however,
it could be much larger. The point of taking the ultralimit is that some
properties of $X$ improve after passing to $X_\om$.

\begin{lem}
$X_\om$ is a geodesic metric space.
\end{lem}
\proof Points $x_\om, y_\om$ in $X_\om$ are represented by sequences
$(x_i), (y_i)$ in $X$. For each $i$ choose a $\frac{1}{i}$-geodesic 
curve $\ga_i$ in $X$ connecting $x_i$ to $y_i$. Then
$$
\ga_\om:=\ulim \ga_i
$$
is a geodesic connecting $x_\om$ to $y_\om$. \qed

\medskip
Similarly, every sequence of
$\frac{1}{i}$-geodesic segments $\ol{y x_i}$ in $X$ satisfying
$$
\ulim d(y, x_i)=\infty,
$$
yields a geodesic ray $\ga_\om$ in $X_\om$ emanating from $y_\om=(y)$.

If $o_i\in X$ is a bounded sequence, then we have a natural (diagonal) isometric 
embedding $X\to X_\om$, given
by the map which sends $x\in X$ to the constant sequence $(x)$.

\begin{lem}\label{ulim}
For every geodesic segment $\ga_\om=\ol{x_\om y_\om}$ in $X_\om$
there exists a sequence of $1/i$-geodesics $\ga_i\subset X_\om$, so
that
$$
\om\hbox{{\rm -lim}}~ \ga_i =\ga_\om.
$$
\end{lem}
\proof Subdivide the segment $\ga_\om$ into $n$ equal subsegments
$$
\ol{z_{\om,j} z_{\om, j+1}}, j=1,...,n,
$$
where $x_\om =z_{\om,1}, y_\om=z_{\om, n+1}$. Then the points
$z_{\om, j}$ are represented by sequences $(z_{k,j})\in X$. It
follows that for $\om$-all $k$, we have
$$
|\sum_{j=1}^n d(z_{k,j}, z_{k, j+1}) - d(x_k, y_k)|< \frac{1}{2i}.
$$
Connect the points $z_{k,j}, z_{k,j+1}$ by $\frac{1}{2i}$-geodesic
segments $\al_{k,j}$. Then the concatenation 
$$
\al_n= \al_{k,1}* ... *\al_{k,n}
$$
is an $\frac{1}{i}$-geodesic connecting $x_k$ and $y_k$, where
$$
x_\om =(x_k), \quad y_\om =(y_k).
$$
It is clear from the construction, that, if given $i$ we choose
sufficiently large $n=n(i)$, then
$$
\ulim \al_{n(i)}=\ga.
$$
Therefore we take $\ga_i:=\al_{n(i)}$. \qed

\subsection{Tripods}

Our next goal is to define {\em tripods} in $X$, which will be our
main technical tool. Suppose that $x, y, z, o$ are points in $X$ and $\mu$ is
an $\eps$-geodesic segment connecting $x$ to $y$, so that $o\in \mu$
and
$$
o\in p_{\mu, \eps}(z).
$$
Then the path $\mu$ is the concatenation $\al\cup \be$, where $\al,
\be$ are $\eps$-geodesics connecting $x, y$ to $o$. Let $\ga$ be an
$\eps$-geodesic connecting $z$ to $o$.

\begin{defn}
1. We refer to $\al\cup \be\cup \ga$ as a {\em tripod} $T$ with the
vertices $x, y, z$, legs $\al, \be, \ga$, and the center $o$.

2.  Suppose that the length of $\al, \be, \ga$ is at least $R$. Then
we refer to the tripod $T$ as $(R,\eps)$-tripod. An $(R,0)$-tripod will be called simply an $R$-tripod. 
\end{defn}

The reader who prefers to work with proper geodesic metric spaces
can safely assume that $\eps=0$ in the above definition and thus $T$
is a geodesic tripod.

\begin{figure}[tbh]
\begin{center}
\input{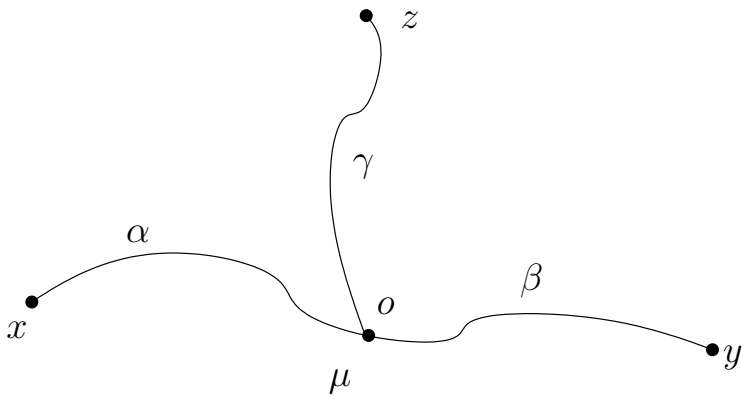}
\end{center}
\caption{\sl A tripod.} \label{tripod.fig}
\end{figure}

\begin{defn}
Let $R\in [0,\infty), \eps\in [0, \infty)$.  We say that a space $X$
is \newline $(R,\eps)$-{\em compressed} if it contains no $(R,\eps)$-tripods. We will
refer to $(R,0)$-compressed spaces as $R$-compressed. A space $X$ which is
not $(R, \eps)$-compressed for any $R<\infty$, $\eps>0$ is called {\em
uncompressed}.
\end{defn}

Therefore $X$ is uncompressed if and only if there exists a sequence of
$(R_i,\eps_i)$-tripods in $X$ with
$$
\lim_i R_i=\infty, \quad \lim_i \eps_i=0.
$$

\subsection{Tripods and ultralimits}

Suppose that $X$ is uncompressed and thus there exists a sequence of
$(R_i,\eps_i)$-tripods $T_i$ in $X$ with
$$
\lim_i R_i=\infty, \quad \lim_i \eps_i=0,
$$
so that the  center of $T_i$ is $o_i$ and the legs are $\al_i,
\be_i, \ga_i$. Then the tripods $T_i$ clearly yield a geodesic
$(\infty, 0)$-tripod $T_\om$ in $(X_\om, o_\om) =\ulim (X, o_i)$. The tripod
$T_\om$ is the union of three geodesic rays $\al_\om, \be_\om,
\ga_\om$ emanating from $o_\om$, so that
$$
o_\om= p_{\mu_\om}(\ga_\om).
$$
Here $\mu_\om =\al_\om \cup \be_\om$. In particular, $X_\om$ is
uncompressed.

Conversely, in view of Lemma \ref{ulim}, we have:

\begin{lem}
\label{pinching} If $X$ is $(R,\eps)$-compressed for $\eps>0$ and
$R<\infty$, then $X_\om$ is $R'$-compressed for every $R'>R$.
\end{lem}
\proof Suppose that $X_\om$ contains an $R'$-tripod $T_\om$.
Then $T_\om$ appears as the ultralimit of $(R'-\frac{1}{i},
\frac{1}{i})$-tripods in $X$. This contradicts the assumption that
$X$ is $(R,\eps)$-compressed. \qed

\medskip
Let $\si: [a,b]\to X$ be a rectifiable curve in $X$ parametererized by its arc-length. 
We let $d_\si$ denote the path metric on $[a,b]$ which is the pull-back of the path
metric on $[a, b]$. By abusing the notation we denote by $d$ the
restriction to $\si$ of the metric $d$. Note that, in general,
$d$ is only a pseudo-metric on $[a, b]$ since $\si$ need not be injective. 
However, if $\si$ is injective then $d$ is a metric.

We repeat this construction with respect to the tripods: Given a
tripod $T\subset X$, define an abstract tripod $T_{mod}$ whose legs
have the same length as the legs of $T$. We have a natural map
$$
\tau: T_{mod}\to X
$$
which sends the legs of $T_{mod}$ to the respective legs of $T$, parameterizing them by the arc-length.
Then $T_{mod}$ has the path metric $d_{mod}$ obtained by pull-back of the
path metric from $X$ via $\tau$. We also have the restriction pseudo-metric
$d$ on $T_{mod}$:
$$
d(A, B)=d(\tau(A), \tau(B)). 
$$
Observe that if $\eps=0$ and $X$ is a tree then
the metrics $d_{mod}$ and $d$ on $T$ agree.

\begin{lem}\label{compare}
$$
d\le d_{mod} \le 3d + 6\eps.
$$
\end{lem}
\proof The inequality $d\le d_{mod}$ is clear. We will prove the
second inequality. If $A, B\in \al\cup \be$ or $A, B\in \ga$ then,
clearly,
$$
d_{mod}(A,B)\le d(A,B)+\eps,
$$
since these curves are $\eps$-geodesics.

Therefore consider the case when $A\in \ga$ and $B\in \be$. Then
$$
D:=d_{mod}(A,B)=t+s,
$$
where $t=d_\al(A,o), s=d_\be(o, B)$.

{\bf Case 1:} $t\ge D/3$. Then, since $o\in \al\cup \be$ is
$\eps$-nearest to $A$, it follows that
$$
D/3\le t\le d(A, o)+\eps \le d(A, B) +2\eps.
$$
Hence
$$
d_{mod}(A,B)=\frac{3D}{3} \le 3(d(A, B) +2\eps)= 3 d(A,B)+ 6\eps,
$$
and the assertion follows in this case.

{\bf Case 2:} $t< D/3$. By the triangle inequality,
$$
D-t=s\le d(o, B)+\eps\le d(o, A)+ d(A, B)+\eps\le t+2\eps+ d(A, B).
$$
Hence
$$
\frac{D}{3}= D- \frac{2}{3}D\le D- 2t \le 2\eps+ d(A, B),
$$
and
$$
d_{mod}(A,B)=\frac{3D}{3}\le 3d(A,B)+ 6\eps. \qed
$$

\section{Topology of configuration spaces of tripods}

We begin with the model tripod $T$ with the legs $\al_i$, $i=1, 2,3$, 
and the center $o$. Consider the configuration space $Z:=T^3\setminus diag$, where $diag$
is the small diagonal
$$
\{(x_1,x_2,x_3)\in T^3: x_1=x_2=x_3\}.
$$
We recall that the homology is taken with the $\Z_2$-coefficients. 

\begin{prop}
$H_1(Z)=0$.
\end{prop}
\proof  $T^3$ is the union of  cubes
$$
Q_{ijk}=\al_i\times \al_j\times \al_k,
$$
where $i, j, k\in \{1, 2,3\}$. Identify each cube $Q_{ijk}$ with the
unit cube in the positive octant in $\R^3$. Then in the cube
$Q_{ijk}$ we choose the equilateral triangle $\si_{ijk}$ given by
the intersection of $Q_{ijk}$ with the hyperplane
$$
x+y+z=1
$$
in $\R^3$. Define the 2-dimensional complex
$$
S:=\bigcup_{ijk} \si_{ijk}.
$$
This complex is homeomorphic to the link of $(o,o,o)$ in $T^3$.
Therefore $Z$ is homotopy-equivalent to
$$
W:=S \setminus (\si_{111}\cup \si_{222} \cup \si_{333}).
$$

Consider the loops $\ga_i:= \partial \si_{iii}$, $i=1, 2, 3$.

\begin{lem}
1. The homology classes $[\ga_i], i=1, 2, 3$ generate $H_1(W)$.

2. $[\ga_1]=[\ga_2]=[\ga_3]$ in  $H_1(W)$.
\end{lem}
\proof 1. We first observe that $S$ is a 2-dimensional spherical
building. Therefore $L$ is homotopy-equivalent to a bouquet of
2-spheres (see \cite[Theorem 2, page 93]{Brown(1989)}), which implies that $H_1(S)=0$. Now the first
assertion follows from the long exact sequence of the pair $(S, W)$.

2. Let us verify that $[\ga_1]=[\ga_2]$. The subcomplex
$$
S_{12}= S \cap (\al_1\cup \al_2)^3
$$
is homeomorphic to the 2-sphere. Therefore $S_{12}\cap W$ is the
annulus bounded by the circles $\ga_1$ and $\ga_2$. Hence
$[\ga_1]=[\ga_2]$. \qed

\begin{lem}
$$
[\ga_1]+[\ga_2]+[\ga_3]=0
$$
in $H_1(W)$.
\end{lem}
\proof Let $B'$ denote the chain
$$
\sum_{\{ijk\}\in A} \si_{ijk},
$$
where $A$ is the set of triples of distinct indices. Let
$$
B'':= \sum_{i=1}^3 (\si_{ii(i+1)}+ \si_{i(i+1)i}+ \si_{(i+1)ii})
$$
where we set $3+1:=1$. We leave it to the reader to verify that
$$
\partial (B'+ B'')=\ga_1+\ga_2+\ga_3
$$
in $C_1(W)$. \qed

By combining these lemmata we obtain the assertion of the
proposition. \qed

\bigskip
{\bf Application to tripods in metric spaces.} 
Consider an $(R,\eps)$-tripod $T$ in a metric space $X$
and its standard parametrization $\tau: T_{mod}\to T$. 

There is an obvious scaling operation 
$$
u\mapsto r\cdot u$$
on the space $(T_{mod}, d_{mod})$ which sends each leg to itself and 
scales all  distances by $r\in [0, \infty)$. It induces the map $T_{mod}^3\to T_{mod}^3$, 
denoted $t\mapsto r\cdot t$, $t\in T_{mod}^3$.


We have the functions
$$
L_{mod}: T^3_{mod}\to K, \quad L: T^3_{mod}\to K,
$$
$$
L_{mod}(x,y,z)= (d_{mod}(x,y), d_{mod}(y,z), d_{mod}(z,x)),
$$
$$
L(x,y,z)= (d(x,y), d(y,z), d(z,x))
$$
computing side-lengths of triangles with respect to the metrics
$d_{mod}$ and $d$.

For $\rho\ge 0$ set
$$
K_\rho:=\{ (a,b,c)\in K: a+b+c>\rho\}.
$$
Define
$$
T^3(\rho):=L^{-1}(K_\rho), \quad T^3_{mod}(\rho):=
L_{mod}^{-1}(K_\rho).
$$
Thus
$$
T^3_{mod}(0)=T^3(0)= T^3\setminus \hbox{~diag}.
$$

\begin{lem}
For every $\rho\ge 0$, the space $T^3_{mod}(\rho)$ is homeomorphic to
$T^3_{mod}(0)$.
\end{lem}
\proof Recall that $S$ is the link of $(o,o,o)$ in $T^3$. Then scaling
determines homeomorphisms
$$
T_{mod}(\rho) \to S \times \R \to T_{mod}(0). \qed
$$

\begin{cor}
For every $\rho\ge 0$, $H_1(T_{mod}(\rho), \Z_2)=0$.
\end{cor}

\begin{cor}\label{zero}
The map induced by inclusion
$$
H_1(T^3(3\rho+18\eps)) \to H_1(T^3(\rho))
$$
is zero.
\end{cor}
\proof Recall that
$$
d\le d_{mod} \le 3 d +6 \eps.
$$
Therefore
$$
T^3(3\rho+ 18 \eps) \subset T^3_{mod}(\rho) \subset T^3(\rho).
$$
Now the assertion follows from the previous  corollary. \qed

\section{Proof of Theorem \ref{B}}

Suppose that $X$ is uncompressed. Then for every $R<\infty, \eps>0$
there exists an $(R,\eps)$-tripod $T$ with the legs $\al, \be, \ga$.
Without loss of generality we may assume that the legs of $T$ have
length $R$. Let $\tau: T_{mod}\to T$ denote the standard map from
the model tripod onto $T$. We will continue with the notation of the previous section. 

Given a map $f: E\to T_{mod}^3$ (or a chain $\si \in C_*(T_{mod}^3)$) 
let $\hat f$ (resp. $\hat\si$) denote the map $L \circ f$ 
from $E$ to $K$ (resp. the chain $L_*(\si)\in C_*(K)$). Similarly, we define $\hat{f}_{mod}$ and $\hat\si_{mod}$ 
using the map $L_{mod}$ instead of $L$.

Every loop $\la: S^1\to T_{mod}^3$, determines the map of the 2-disk
$$
\La: D^2 \to T_{mod}^3,
$$
given by
$$
\La(r, \theta)= r\cdot \la(\theta)
$$
where we are using the polar coordinates $(r, \theta)$ on the unit disk $D^2$. 
Triangulating both $S^1$ and $D^2$ and assigning the coefficient $1\in \Z_2$ to each simplex,  
we regard both $\la$ and $\La$ as singular chains in $C_*(T_{mod}^3)$.

We let $a, b, c$ denote the coordinates on the space $\R^3$ containing the cone $K$. 
Let $\kappa=(a_0, b_0, c_0)$ be a $\del$-nondegenerate point in the
interior of $K$ for some $\del>0$; set $r:= a_0+b_0+c_0$. 

Suppose that there exists a loop $\la$ in $T_{mod}^3$ such that:

1. $\hat\la(\theta)$ is $\eps$-degenerate for each $\theta$. Moreover, each triangle $\la(\theta)$ 
is either contained in $\al_{mod}\cup \be_{mod}$ or has only two distinct vertrices. 

In particular, the image of $\hat\la$ is contained in 
$$
K \setminus \R_+\cdot \kappa. 
$$

2. The image of $\hat\la$ is contained in $K_\rho$, where $\rho= 3r+18\eps$.

3. The homology class $[\hat\la]$ is nontrivial in 
$H_1(K \setminus \R_+\cdot \kappa)$.

\begin{figure}[tbh]
\begin{center}
\input{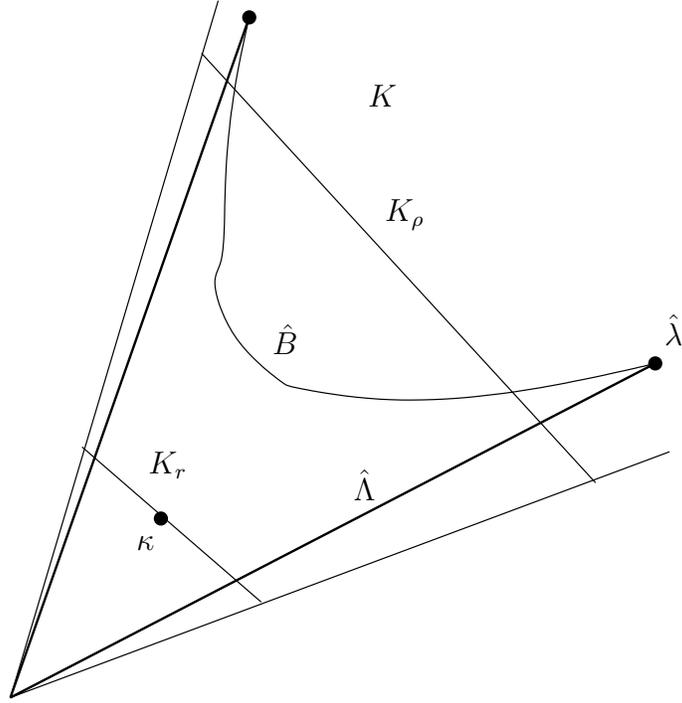}
\end{center}
\caption{\sl Chains $\hat\La$ and $\hat{B}$.} \label{loop.fig}
\end{figure}

\begin{lem}
\label{belongs}
If there exists a loop $\la$ satisfying the assumptions 1---3, and $\eps<\del/2$, then
$\kappa$ belongs to $K_3(X)$.
\end{lem}
\proof We have the 2-chains
$$
\hat\La, \hat\La_{mod} \in C_2(K \setminus \kappa),
$$
with
$$
\hat\la =\partial \hat \La, \hat\la_{mod}=\partial \hat \La_{mod} \in C_1(K_\rho).
$$
Note that the support of $\hat\la_{mod}$ is contained in $\partial K$ and 
the 2-chain $\hat\La_{mod}$ is obtained by coning off $\hat\la_{mod}$ from the origin. 
Then, by Assumption 1, for every $w\in D^2$:

i. Either $d( \hat \La(w), \hat \La_{mod}(w))\le \eps$. 

ii. Or $\hat \La(w), \hat \La_{mod}(w)$ belong to the common ray in $\partial K$. 

Since $d(\kappa, \partial K)> \del\ge 2\eps$, it follows that the straight-line homotopy $H_t$ 
between the maps  
$$
\hat \La, \hat \La_{mod}: D^2\to K$$
misses $\kappa$. Since $K_\rho$ is convex, $H_t(S^1)\subset K_\rho$ for each $t\in [0,1]$, and 
we obtain 
$$
[\hat\La_{mod}]=[\hat\La]\in H_2(K\setminus \kappa, K_\rho). 
$$
Assumptions 2 and 3 imply that the relative homology class
$$
[\hat\La_{mod}]\in H_2(K\setminus \kappa, K_\rho)
$$
is nontrivial. Hence
$$
[\hat\La]\in H_2(K\setminus \kappa, K_\rho)
$$
is nontrivial as well. Since $\rho=3r+18\eps$, according to Corollary \ref{zero}, $\la$ bounds a 2-chain 
$$
{\mathrm B}\in C_2(T^3(r)). 
$$
Set $\Si:= {\mathrm B}+\La$. 
Then the absolute class 
$$
[\hat\Si]=[\hat\La +\hat {\mathrm B}]\in H_2(K \setminus \kappa)$$
is also nontrivial.  Since $T_{mod}^3$ is contractible, 
there exists a 3-chain $\Ga\in C_3(T_{mod}^3)$ such that
$$
\partial \Ga=\Si.
$$
Therefore the support of $\hat\Ga$ contains the point $\kappa$. Since the
map
$$
L: T^3\to K
$$
is the composition of the continuous map $\tau^3: T^3\to X^3$ with the
continuous map $L: X^3\to K$, it follows that $\kappa$ belongs to
the image of the map $L: X^3\to K$ and hence $\kappa\in K_3(X)$.
\qed

Our goal therefore is to construct a loop $\la$,
satisfying Assumptions 1---3.

\medskip
Let $T\subset X$ be an $(R,\eps)$-tripod with the legs $\al, \be,
\ga$ of the length $R$, where $\eps\le \del/2$.  We let $\tau: T_{mod}\to T$ denote the
standard parametrization of $T$. Let $x, y, z, o$ denote the
vertices and the center of $T_{mod}$. We let $\al_{mod}(s),
\be_{mod}(s), \ga_{mod}(s): [0, R]\to T_{mod}$ denote the arc-length
parameterizations of the legs of $T_{mod}$, so that
$\al(R)=\be(R)=\ga(R)=o$. 

We will describe the loop $\la$ as the concatenation of seven paths 
$$
p_i(s)=(x_1(s), x_2(s), x_3(s)), i=1,...,7.
$$
We let $a=d(x_2, x_3), b=d(x_3,x_1), c=d(x_1, x_2)$.

1. $p_1(s)$ is the path starting at $(x, x, o)$ and ending at $(o, x, o)$, given by
$$
p_1(s)= (\al_{mod}(s), x, o).
$$
Note that for $p_1(0)$ and $p_1(R)$ we have $c=0$ and $b=0$ respectively.

2. $p_2(s)$ is the path starting at $(o, x, o)$ and ending at $(y, x, o)$, given by
$$
p_2(s)= (\bar\be_{mod}(s), x, o).
$$

3. $p_3(s)$ is the path starting at $(y, x, o)$ and ending at $(y, o, o)$, given by
$$
p_3(s)= (y, \al_{mod}(s), o).
$$
Note that for $p_3(R)$ we have $a=0$.

4. $p_4(s)$ is the path starting at $(y, o, o)$ and ending at $(y, y, o)$, given by
$$
p_4(s)= (y, \bar\be_{mod}(s), o).
$$

Note that for $p_4(R)$ we have $c=0$. Moreover, if $\al*\bar\be$ is a geodesic, then
$$
d(\tau(x), \tau(o))=d(\tau(y), \tau(o)) \Rightarrow \hat{p}_4(R)=\hat{p}_1(0)
$$
and therefore $\hat{p}_1*... * \hat{p}_4$ is a loop.

5. $p_5(s)$ is the path starting at $(y, y, o)$ and ending at $(y, y, z)$ given by
$$
(y, y, \bar\ga_{mod}(s)).
$$

6. $p_6(s)$ is the path starting at $(y, y, z)$ and ending at $(x, x, z)$ given by
$$
(\be_{mod} * \bar\al_{mod}, \be_{mod} * \bar\al_{mod}, z).
$$

7. $p_7(s)$ is the path starting at $(x, x, z)$ and ending at $(x, x, o)$ given by
$$
(x, x, \ga_{mod}(s)).
$$

Thus
$$
\la:= p_1* ... * p_7
$$
is a loop.

Since $\al*\be$ and $\ga$ are $\eps$-geodesics in $X$, each path
$p_i(s)$ determines a family of $\eps$-degenerate triangles in
$(T_{mod}, d)$. It is clear that Assumption 1 is satisfied.

\begin{figure}[tbh]
\begin{center}
\input{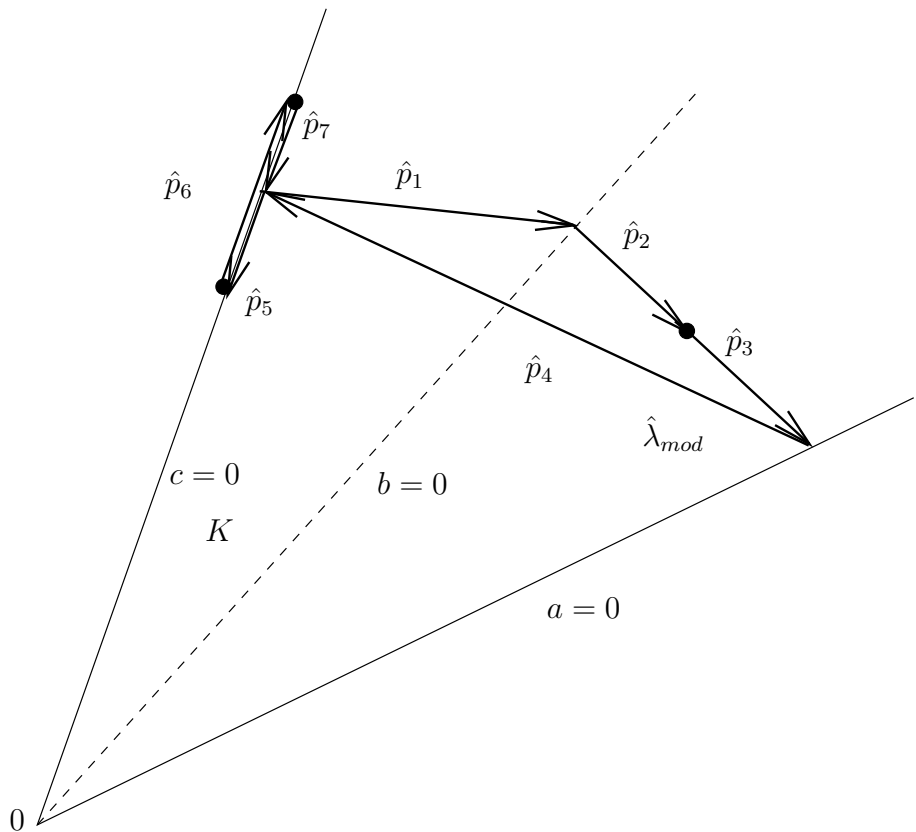}
\end{center}
\caption{\sl The loop $\hat\la_{mod}$.} \label{f1.fig}
\end{figure}

The class $[\hat\la_{mod}]$ is clearly nontrivial in $H_1(\partial K\setminus 0)$. 
See Figure \ref{f1.fig}. Therefore, since $\eps\le \del/2$, 
$$
[\hat\la]=[\hat\la_{mod}]\in H_1(K \setminus \R_+\cdot \kappa)\setminus \{0\},
$$
see the proof of Lemma \ref{belongs}. Thus Assumption 2 holds.

\begin{lem}
The image of $\hat\la$ is contained in the closure of $K_{\rho'}$, where
$$
\rho'= \frac{2}{3}R -4\eps.
$$
\end{lem}
\proof We have to verify that for each $i=1,...,6$ and every $s\in
[0, R]$, the perimeter (with respect to the metric $d$) of each
triangle $p_i(s)\in T_{mod}^3$ is at least $\rho'$. These
inequalities follow directly from Lemma \ref{compare} and the
description of the paths $p_i$. \qed

\medskip
Therefore, if we take 
$$
R> \frac{9}{2} r - 33\eps
$$
then the image of $\hat\la$ is contained in 
$$
K_{3r+18\eps}
$$
and Assumption 3 is satisfied. Theorem \ref{A} follows. \qed

\section{Quasi-isometric characterization of compressed spaces}

The goal of this section is to prove Theorem \ref{B}. Suppose that
$X$ is compressed. The proof is easier if $X$ is a proper
geodesic metric space, in which case there is no need considering
the ultralimits. Therefore, we recommend the reader uncomfortable
with this technique to assume that $X$ is a proper geodesic metric
space.

Pick a base-point $o\in X$, a nonprincipal ultrafilter $\om$ and
consider the ultralimit
$$
X_\om =\ulim (X, o)
$$
of the constant sequence of pointed metric spaces. If $X$ is a
proper geodesic metric space then, of course, $X_\om =X$.
In view of Lemma \ref{pinching}, the space $X_\om$ is $R$-compressed
for some $R$. 

Assume that $X$ is unbounded. Then $X$ contains a 
sequence of $1/i$-geodesic paths $\ga_i=\ol{o x_i}$ with
$$
\ulim d(o, x_i)=\infty,
$$
which yields a geodesic ray $\rho_1$ in $X_\om$ emanating from the point $o_\om$.

\begin{lem}\label{ray}
Let $\rho$ be a geodesic ray in $X_\om$ emanating from the point $O$. 
Then the neighborhood $E=N_R(\rho)$ is an end $E(\rho)$ of $X_\om$. 
\end{lem}
\proof Suppose that $\al$ is a
path in $X_\om\setminus B_{2R}(O)$ connecting a point $y\in X_\om\setminus
E$ to a point $x\in E$. Then there exists a point $z\in \al$ such
that $d(z, \rho)=R$. Since $X_\om$ contains no $R$-tripods,
$$
d(p_{\rho}(z), O)<R.
$$
Therefore $d(z, O)<2R$. Contradiction. \qed

\medskip
Set $E_1:=E(\rho_1)$. If the image of the natural embedding $\iota: X\to X_\om$ is contained in a finite 
metric neighborhood of $\rho_1$, then we are done, as $X$ is quasi-isometric to $\R_+$.  
Otherwise, there exists a sequence $y_n\in X$ such that:
$$
\ulim d(\iota(y_n), \rho_1)=\infty.
$$
Consider the $\frac{1}{n}$-geodesic paths $\al_n\in P(o, y_n)$. 
The sequence $(\al_n)$ determines a geodesic ray $\rho_2\subset X_\om$ emanating from $o_\om$. 
Then there exists $s\ge 4R$ such that 
$$
d(\al_n(s), \ga_i)\ge 2R
$$
for $\om$-all $n$ and $\om$-all $i$. Therefore, for $t\ge s$, $\rho_2(t)\notin E(\rho_1)$. 
By applying Lemma \ref{ray} to $\rho_2$ we conclude that $X_\om$ has an end 
$E_2=E(\rho_2)=N_R(\rho_2)$. Since $E_1, E_2$ are distinct ends of $X_\om$, 
$E_1\cap E_2$ is a bounded subset. Let $D$ denote the diameter of this intersection. 

\begin{lem}\label{qg}
1. For every pair of points $x_i=\rho_i(t_i)$, $i=1, 2$, we have
$$
\ol{x_1 x_2}\subset N_{D/2+ 2R}(\rho_1\cup \rho_2).
$$

2. $\rho_1\cup \rho_2$ is a quasi-geodesic.
\end{lem}
\proof Consider the points $x_i$ as in Part 1. Our goal is to get a
lower bound on $d(x_1, x_2)$. A geodesic segment $\ol{x_1 x_2}$
has to pass through the ball $B(o_\om, 2R), i=1, 2$, since this ball
separates the ends $E_1, E_2$. Let $y_i\in \ol{x_1 x_2}\cap B(o_\om,
2R)$ be such that
$$
\ol{x_i y_i}\subset E_i, i=1, 2.
$$
Then
$$
d(y_1, y_2)\le D+ 4R, 
$$
$$
d(x_i, y_i)\ge t_i - 2R,
$$
and 
$$
\ol{x_i y_i}\subset N_{R}(\rho_i), \quad i=1, 2.
$$
This implies the first assertion of Lemma. Moreover,
$$
d(x_1, x_2)\ge d(x_1, y_1)+ d(x_2, y_2)\ge t_1+t_2- 4R =d(x_1,
x_2)-4R.
$$
Therefore $\rho_1\cup \rho_2$ is a $(1, 4R)$--quasi-geodesic.
\qed

If $\iota(X)$ is contained in a finite metric neighborhood of $\rho_1\cup \rho_2$, then, 
by Lemma \ref{qg}, $X$ is quasi-isometric to $\R$. Otherwise, there exists a 
sequence $z_k\in X$ such that 
$$
\ulim d(\iota(z_k), \rho_1\cup \rho_2)=\infty.
$$
By repeating the construction of the ray $\rho_2$, we obtain a geodesic ray $\rho_3\subset X_\om$ 
emanating from the point $o_\om$, so that $\rho_3$ is not contained in a finite metric neighborhood of  
$\rho_1\cup \rho_2$. For every $t_3$,
the nearest-point projection of $\rho_3(t_3)$ to
$$
N_{D/2+ 2R}(\rho_1\cup \rho_2)$$
is contained in
$$
B_{2R}(o_\om).
$$
Therefore, in view of Lemma \ref{qg}, for every pair of points
$\rho_i(t_i)$ as in that lemma, the nearest-point projection of
$\rho_3(t_3)$ to $\ol{\rho_1(t_1)\rho_2(t_2)}$ is contained in
$$
B_{4R+D}(o_\om).
$$
Hence, for sufficiently large $t_1, t_2, t_3$, the points
$\rho_i(t_i)$, $i=1, 2, 3$ are vertices of an $R$-tripod in $X$. 
This contradicts the assumption that $X_\om$ is $R$-compressed. 

Therefore $X$ is either bounded, or is quasi-isometric to a $\R_+$ or to $\R$. \qed

\section{Examples}
\label{examples}

\begin{thm}\label{example1}
There exist an (incomplete) 2-dimensional Riemannian manifold $M$ quasi-isometric to $\R$, so that:

a.  $K_3(M)$ does not contain $\partial K_3(\R^2)$.

b. For the Riemannian product $M^2=M\times M$, $K_3(M^2)$ does not
contain $\partial K_3(\R^2)$ either.

Moreover, there exists $D<\infty$ such that for every degenerate triangle in $M$ and $M^2$, at 
least one side is $\le D$. 
\end{thm}
\proof a. We start with the open concentric annulus $A\subset \R^2$, which has
the inner radius $R_1>0$ and the outer radius $R_2<\infty$. We give
$A$ the flat Riemannian metric induced from $\R^2$. Let $M$ be
the universal cover of $A$, with the pull-back Riemannian
metric. Since $M$ admits a properly discontinuous isometric action of $\Z$ with the 
quotient of finite diameter, it follows that 
$M$ is quasi-isometric to $\R$. The metric completion  $\bar{M}$ of $M$
is diffeomorphic to the closed bi-infinite flat  strip. 
Let $\partial_1 M$ denote the component of the boundary of
$\bar{M}$ which covers the inner boundary of $A$ under the map of
metric completions
$$
\bar{M}\to \bar{A}.
$$
As a metric space, $\bar{M}$ is $CAT(0)$, therefore it contains a
unique geodesic between any pair of points. However, for any pair of
points $x, y\in M$, the geodesic $\ga=\ol{xy}\subset \bar{M}$ is the union of subsegments
$$
\ga_1\cup \ga_2 \cup \ga_3
$$
where $\ga_1, \ga_3\subset M$, $\ga_2\subset \partial_1 M$, and the
lengths of $\ga_1, \ga_3$ are at most $D_0=\sqrt{R_2^2- R_1^2}$.

Hence, for every degenerate triangle $(x,y,z)$ in $M$, at least one side is $\le D_0$.

\begin{figure}[tbh]
\begin{center}
\input{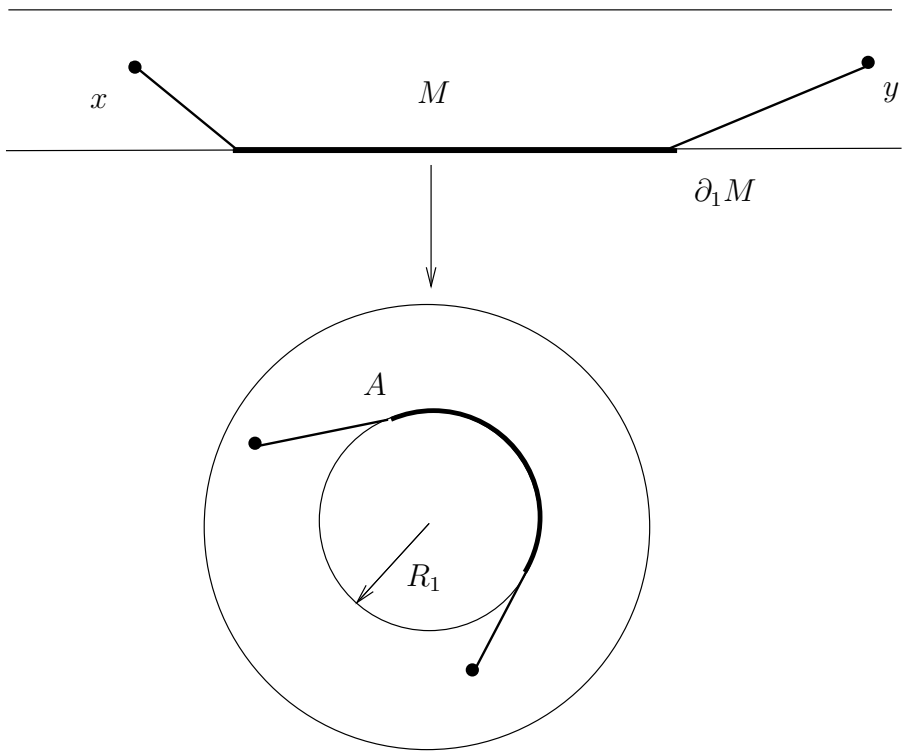}
\end{center}
\caption{\sl Geodesics in $\bar{M}$.} \label{f2.fig}
\end{figure}

\medskip
b. We observe that the metric completion of $M^2$ is $\bar{M}\times
\bar{M}$; in particular, it is again a $CAT(0)$ space. Therefore it
has a unique geodesic between any pair of points. Moreover,
geodesics in $\bar{M}\times \bar{M}$ are of the form
$$
(\ga_1(t), \ga_2(t))
$$
where $\ga_i$, $i=1,2$ are geodesics in $\bar{M}$. Hence for
every geodesic segment $\ga\subset \bar{M}\times \bar{M}$, the complement
$\ga\setminus \partial \bar{M}^2$ is the union of two subsegments of
length $\le \sqrt{2}D_0$ each. Therefore for every degenerate triangle in $M^2$, 
at least one side is $\le \sqrt{2}D_0$. \qed

\begin{rem}
The manifold $M^2$ is, of course, quasi-isometric to $\R^2$.
\end{rem}

Our second example is a graph-theoretic analogue of the
Riemannian manifold $M$.

\begin{thm}\label{example2}
There exists a complete path metric space $X$ (a metric
graph)  quasi-isometric to $\R$ so that:

a. $K_3(X)$ does not contain $\partial K_3(\R^2)$.

b. $K_3(X^2)$ does not contain $\partial K_3(\R^2)$.

Moreover, there exists $D<\infty$ such that for every degenerate triangle in $X$ and $X^2$, at 
least one side is $\le D$. 
\end{thm}
\proof  a. We start with the disjoint union of oriented circles $\al_i$ of the
length $1+\frac{1}{i}$, $i\in I=\N \setminus \{2\}$. We regard each $\al_i$ as a path
metric space. For each $i$ pick a point $o_i\in \al_i$  and its
antipodal point $b_i\in \al_i$. We let $\al_i^+$ be the positively
oriented arc of $\al_i$ connecting $o_i$ to $b_i$. Let $\al_i^-$ be
the complementary arc.

Consider the bouquet $Z$ of $\al_i$'s by gluing them all at the points
$o_i$. Let $o\in Z$ be the image of the points $o_i$. Next, for every pair 
$i, j\in I$ attach to $Z$ the oriented arc $\be_{ij}$ of the length
$$
\frac{1}{2} + \frac{1}{4}( \frac{1}{i} +  \frac{1}{j})
$$
connecting $b_i$ and $b_j$ and oriented from $b_i$ to $b_j$ if
$i<j$. Let $Y$ denote the resulting graph. We give $Y$ the path
metric. Then $Y$ is a complete metric space, since it is a metric
graph where the length of every edge is at least $1/2>0$. 
Note also that the length of every edge in $Y$ is at most $1$. 

\begin{figure}[tbh]
\begin{center}
\input{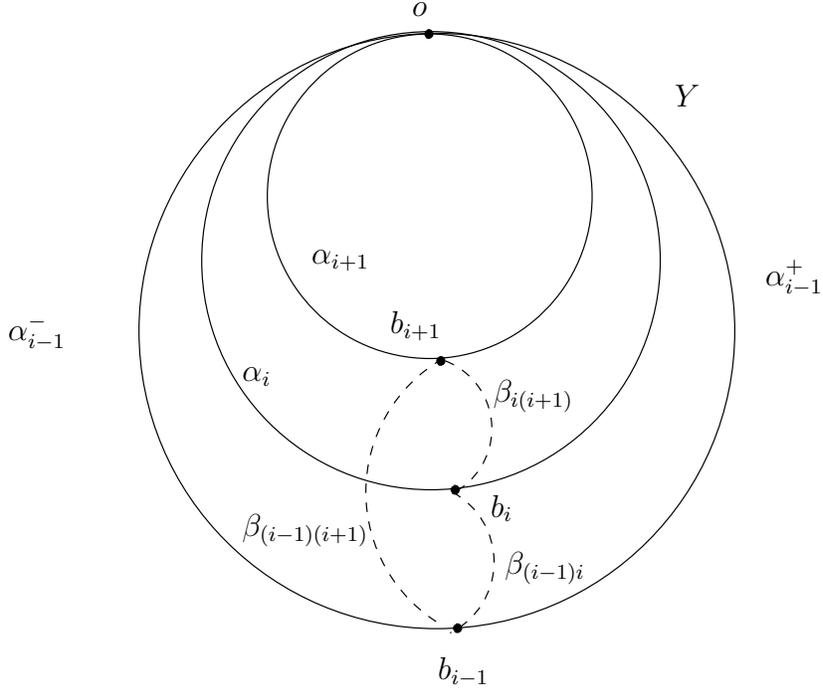}
\end{center}
\caption{\sl The metric space $Y$.} \label{f4.fig}
\end{figure}

The space $X$ is the infinite cyclic regular cover over $Y$ defined as
follows. Take the maximal subtree
$$
T=\bigcup_{i\in I} \al_i^+\subset Y.
$$
Every oriented edge of $Y\setminus T$ determines a free generator of
$G=\pi_1(Y, o)$. Define the homomorphism $\rho: G\to \Z$ by sending 
every free generator to $1$. Then the covering $X\to Y$ is
associated with the kernel of $\rho$. (This covering exists since
$Y$ is locally contractible.)

We lift the path metric from $Y$ to $X$, thereby making $X$ a
complete metric graph. We label vertices and edges of $X$ as
follows.

1. Vertices $a_n$ which project to $o$. The cyclic group $\Z$ acts
simply transitively on the set of these vertices thereby giving them the 
indices $n\in \Z$.

2. The edges $\al_i^{\pm}$ lift to the edges $\al_{in}^{+},
\al_{in}^{-}$ incident to the vertices $a_n$ and $a_{n+1}$
respectively.

3. The intersection $\al_{in}^+\cap \al_{i(n+1)}^-$ is the vertex
$b_{in}$ which projects to the vertex $b_i\in \al_i$.

4. The edge $\be_{ijn}$ connecting $b_{in}$ to $b_{j(n+1)}$ which projects 
to the edge $\be_{ij}\subset Y$.

\begin{figure}[tbh]
\begin{center}
\input{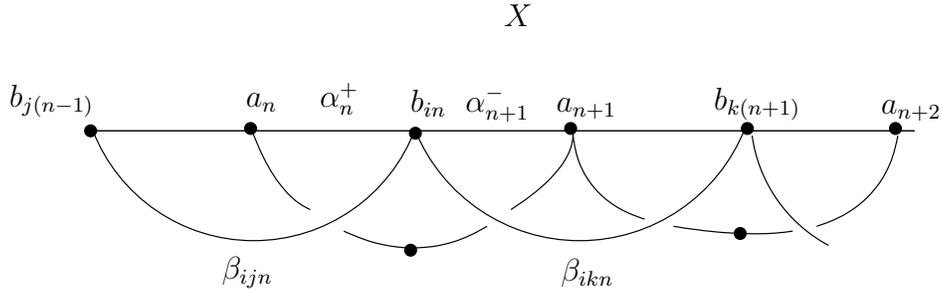}
\end{center}
\caption{\sl The metric space $X$.} \label{f5.fig}
\end{figure}

\begin{lem}
\label{nodeg}
$X$ contains no degenerate triangles $(x, y, v)$, so that $v$ is a
vertex,
$$
d(x,v)+d(v,y)=d(x,y)
$$
and $\min( d(x, v), d(v,y))>2$.
\end{lem}
\proof Suppose that such degenerate triangles exist. 

{\bf Case 1:} $v=b_{in}$. Since the triangle $(x,y,v)$ is
degenerate, for all sufficiently small $\eps>0$ there exist
$\eps$-geodesics $\si$ connecting $x$ to $y$ and passing through
$v$.

Since $d(x, v), d(v, y)>2$, it follows that for sufficiently small
$\eps>0$, $\si=\si(\eps)$ also passes through $b_{j(n-1)}$ and
$b_{k(n+1)}$ for some $j, k$ depending on $\si$. We will assume that
as $\eps\to 0$, both $j$ and $k$ diverge to infinity, leaving the
other cases to the reader.

Therefore
$$
d(x, v)= \lim_{j\to \infty} (d(x, b_{j(n-1)}) +  d(b_{j(n-1)}, v)),
$$
$$
d(v, y)= \lim_{k\to \infty} (d(y, b_{k(n+1)}) +  d(b_{k(n+1)}, v)).
$$
Then
$$
\lim_{j\to \infty} d(b_{j(n-1)}, v) + \lim_{k\to \infty}
d(b_{k(n+1)}, v) = 1+ \frac{1}{2i}.
$$
On the other hand, clearly,
$$
\lim_{j, k\to \infty} d(b_{j(n-1)}, b_{k(n+1)}) = 1.
$$
Hence
$$
d(x, y)= \lim_{j\to \infty} d(x, b_{j(n-1)}) + \lim_{k\to \infty}
d(y, b_{k(n+1)}) +1 < d(x, v) + d(v, y).
$$
Contradiction.

{\bf Case 2:} $v= a_n$. Since the triangle $(x, y, v)$ is
degenerate, for all sufficiently small $\eps>0$ there exist
$\eps$-geodesics $\si$ connecting $x$ to $y$ and passing through
$v$. Then for sufficiently small $\eps>0$, every $\si$ also passes
through $b_{j(n-1)}$ and $b_{kn}$ for some $j, k$ depending on
$\si$. However, since $j, k\ge 2$, 
$$
d(b_{j(n-1)}, b_{kn})= \frac{1}{2}+ \frac{1}{4j} + \frac{1}{4i} 
\le \frac{3}{4} < 1=\inf_{j, k} (d(b_{j(n-1)},
v)+ d(v, b_{kn})).
$$
Therefore $d(x, y)<d(x, v)+ d(v, y)$. Contradiction. \qed

\begin{cor}
$X$ contains no degenerate triangles $(x, y, z)$, such that
$$
d(x,z)+d(z,y)=d(x,y)
$$
and $\min( d(x, z), d(z,y))\ge 3$.
\end{cor}
\proof Suppose that such a degenerate triangle exists. We can assume
that $z$ is not a vertex. The point $z$ belongs to an edge $e\subset
X$. Since $\length(e)\le 1$, for one of the vertices $v$ of $e$
$$
d(z,v)\le 1/2. 
$$
Since the triangle $(x,y,z)$ is degenerate,
for all $\eps$-geo\-de\-sics $\si\in P(x,z)$, $\eta\in P(z,y)$ we have:
$$
e\subset \si\cup \eta,
$$
provided that $\eps>0$ is sufficiently small.  Therefore the
triangle $(x,y,v)$ is also degenerate. Clearly,
$$
\min (d(x, v), d(y, v))\ge \min (d(x, z), d(y, z)) -1/2\ge 2.5.
$$
This contradicts Lemma \ref{nodeg}.  \qed

Hence part (a) of Theorem \ref{example2} follows.

b. We consider $X^2=X\times X$ with the product metric
$$
d^2((x_1, y_1), (x_2, y_2))= d^2(x_1, x_2)+ d^2(y_1, y_2).
$$
Then $X^2$ is a complete path-metric space. Every degenerate triangle in $X^2$ 
projects to degenerate triangles in both factors. It therefore follows from part (a) that $X$ 
contains no degenerate triangles with all sides $\ge 18$. We leave the detals to the reader. \qed

\section{Exceptional cases}
\label{exceptional}

\begin{thm}
\label{L=1}
Suppose that $X$ is a path metric space quasi-isometric to a metric space $X'$, which is either
$\R$ or $\R_+$. Then there exists a $(1, A)$-quasi-isometry $X'\to X$.
\end{thm}
\proof We first consider the case $X'=\R$.
The proof is simpler if $X$ is proper, therefore we sketch it first under this assumption.
Since $X$ is quasi-isometric to $\R$, it is 2-ended with the ends $E_+, E_-$.
Pick two divergent sequences $x_i\in E_+, y_i\in E_-$. Then there exists a compact subset $C\subset X$ so that 
all geodesic segments $\ga_i:=\ol{x_i y_i}$ intersect $C$. 
It then follows from the Arcela-Ascoli theorem that the sequence of segments $\ga_i$ subconverges
to a complete geodesic $\ga\subset X$. Since $X$ is quasi-isometric to $\R$, there exists
$R<\infty$ such that $X=N_R(\ga)$. We define the $(1, R)$-quasi-isometry $f: \ga\to X$ to be the
identity (isometric) embedding. 

\medskip
We now give a proof in the general case. 
Pick a non-principal ultrafilter $\om$ on $\N$ and a base-point
$o\in X$. Define $X_\om$ as the $\om$-limit of $(X, o)$. The quasi-isometry $f: \R\to X$ yields 
a quasi-isometry $f_\om: \R=\R_\om \to X_\om$. Therefore $X_\om$ is also quasi-isometric to $\R$. 

We have the natural isometric embedding $\iota: X\to X_\om$. As above, let $E_+, E_-$ denote the ends
of $X$ and choose divergent sequences $x_i\in E_+, y_i\in E_-$. Let $\ga_i$ denote an
$\frac{1}{i}$-geodesic  segment in $X$ connecting $x_i$ to $y_i$. Then each $\ga_i$ intersects
a bounded subset $B\subset X$. Therefore, by taking the ultralimit of $\ga_i$'s,  we obtain
a complete geodesic $\ga\subset X_\om$. Since $X_\om$ is quasi-isometric to $\R$, the 
embedding $\eta: \gamma\to X_\om$ is a quasi-isometry. Hence $X_\om= N_R(\ga)$ for
some $R<\infty$. 

\medskip
For the same reason, 
$$
X_\om=N_{D}(\iota(X))$$
for some $D<\infty$. Therefore the isometric embeddings 
$$
\eta: \ga \to X_\om, \iota: X\to X_\om$$ 
are $(1,R)$ and $(1,D)$-quasi-isometries respectively. 
By composing $\eta$ with the quasi-inverse to $\iota$, 
we obtain a $(1, R+3D)$-quasi-isometry $\R\to X$. 

The case when $X$ is quasi-isometric to $\R_+$ can be treated as follows. Pick a point $o\in X$ 
and glue two copies of $X$ at $o$. Let $Y$ be the resulting path metric space. It is easy to see that $Y$ 
is quasi-isometric to $\R$ and the inclusion $X\to Y$ is an isometric embedding. 
Therefore, there  exists a $(1, A)$-quasi-isometry $h: Y\to \R$ and  
the restriction of $h$ to $X$ yields the $(1, A)$-quasi-isometry from $X$ to the half-line. \qed

\begin{cor}
\label{C}
Suppose that $X$ is a path metric space quasi-isometric to
$\R$ or $\R_+$. Then $K_3(X)$ is contained in the $D$-neighborhood of $\partial K$ for some $D<\infty$.
In particular, $K_3(X)$ does not contain the interior of $K=K_3(\R^2)$.
\end{cor}
\proof Suppose that $f: X\to X'$ is an $(L,A)$-quasi-isometry, where $X'$ is either $\R$ or $\R_+$.
According to Theorem \ref{L=1}, we can assume that $L=1$.
For every triple of points $x, y, z\in X$, after relabelling, we obtain
$$
d(x, y)+d(y, z)\le  d(x, z)+ D,
$$
where $D=3A$. Then every triangle in $X$ is $D$-degenerate. Hence
$K_3(X)$ is contained in the $D$-neighborhood of $\partial K$. \qed

\begin{rem}
One can construct a metric space $X$ quasi-isometric to $\R$ such that $K_3(X)=K$. Moreover, $X$ 
is isometric to a curve in $\R^2$ (with the metric obtained by the restriction of the metric 
on $\R^2$). Of course, the metric on $X$ is not a path metric. 
\end{rem}

\begin{cor}
Suppose that $X$ is a path metric space. Then the following are equivalent:

1. $K_3(X)$ contains the interior of $K=K_3(\R^2)$.

2. $X$ is not quasi-isometric to the point, $\R_+$ and $\R$. 

3. $X$ is uncompressed.

\end{cor}
\proof 1$\Rightarrow$2 by Corollary \ref{C}. 2$\Rightarrow$3 by Theorem \ref{B}. 3$\Rightarrow$1 by Theorem \ref{A}. 
\qed

\begin{rem}
The above corollary remains valid under the following 
assumption on the metric on $X$, which is weaker than 
being a path metric: 

For every pair of points $x, y\in X$ and every $\eps>0$, there exists a 
$(1,\eps)$-quasi-geodesic path $\al\in P(x,y)$. 
\end{rem}

\bibliography{$HOME/BIB/lit}
\bibliographystyle{siam}

\end{document}